\documentclass[12pt,a4paper,english]{article}
\usepackage{float}
\usepackage[slantedGreek]{mathpazo}
\usepackage[pdftex]{graphicx}
\usepackage{amsfonts}
\usepackage{setspace}
\usepackage{amssymb}
\usepackage{amsthm}
\usepackage{graphicx}
\usepackage{amsmath,scalerel}%
\usepackage{subcaption}
\usepackage{color}
\usepackage[english]{babel}
 
\usepackage{natbib}
\usepackage{epstopdf}

\begin{document}

\title{\vspace{-0.5cm} Bayesian Hypothesis Testing: Redux}

\author{Hedibert F. Lopes and Nicholas G. Polson\\{\em Insper and Chicago Booth}\footnote{
Hedibert Lopes is Professor of Statistics and Econometrics at Insper Institute of Education and Research ({\tt hedibertfl@insper.edu.br}) and
Nicholas Polson is Professor of Econometrics and Statistics
at University of Chicago Booth School of Business ({\tt ngp@chicagobooth.edu}).}}

\date{First draft: March 2018\\
This draft: August 2018}

\maketitle

\begin{abstract}
\noindent Bayesian hypothesis testing is re-examined from the perspective of an \emph{a priori} assessment of the test statistic distribution under the alternative.  By assessing the distribution of an observable test statistic, rather than prior parameter values, we provide a practical default Bayes factor which is straightforward to interpret.  To illustrate our methodology, we provide examples where evidence for a Bayesian strikingly supports the null, but leads to rejection under a classical test. Finally, we conclude with directions for future research.
\end{abstract}

\vspace{0.2in} {\bf Keywords:} Bayesian, Hypothesis testing, Bayes factor, $p$-value, Test statistic, Multiple comparisons.

\section{Introduction}
Bayesians and Classicists are sharply divided on the question of hypothesis testing. 
Hypothesis testing is a cousin to model selection and in a world of high dimensional selection problems, hypothesis testing is as relevant today as it ever has been.  We contrast these two approaches, by re-examining the construction of a hypothesis test, 
motivated by the seminal paper of Edwards, Lindman and Savage (1996) (hereafter ELS) who provide the following
contrast:
\begin{quote}
{\it 
We now show informally,
as much as possible from a
classical point of view, how evidence
that leads to classical rejection of a
null hypothesis at the 0.05 level can
favor that null hypothesis. The loose
and intuitive argument can easily be
made precise. Consider a two-tailed t test
with many degrees of freedom. If a
true null hypothesis is being tested,
$t$ will exceed $1.96$ with probability
2.5\% and will exceed 2.58 with
probability 0.5\%. (Of course, 1.96 and
2.58 are the 5\% and 1\% two-tailed
significance levels; the other 2.5\% and
0.5\% refer to the possibility that t may
be smaller than -1.96 or -2.58.)
So on 2\% of all occasions when true
null hypotheses are being tested, $t$
will lie between 1.96 and 2.58. How
often will $t$ lie in that interval when
the null hypothesis is false? That
depends on what alternatives to the
null hypothesis are to be considered.
Frequently, given that the null hypothesis
is false, all values of $t$ between,
say, $-20$ and $+20$ are about
equally likely for you. Thus, when
the null hypothesis is false, t may well
fall in the range from 1.96 to 2.58 with
at most the probability (2.58 - 1.96)/
[+20 - (-20)] = 1.55\%. In such
a case, since 1.55 is less than 2 the
occurrence of t in that interval speaks
mildly for, not vigorously against, the
truth of the null hypothesis.
This argument, like almost all the
following discussion of null hypothesis
testing, hinges on assumptions about
the prior distribution under the alternative
hypothesis. The classical statistician
usually neglects that distribution
in fact, denies its existence.
He considers how unlikely a t as far
from 0 as 1.96 is if the null hypothesis
is true, but he does not consider that
a $t$ as close to 0 as 1.96 may be even
less likely if the null hypothesis is
false.}
\end{quote}

In terms of a decision rule\footnote{Here $\Omega(A)$ is the prior odds of the null.  $\Omega(A|D)$ is the 
posterior odds given datum $D$, and $L(A;D)$ is the likelihood ratio (a.k.a. Bayes factor, BF).}, ELS go on to say:

\begin{quote}
{\it 
If you need not make your guess until after you have examined a datum $D$, you will prefer to guess $A$ if and only if $ \Omega(A|D) $ exceeds
$J/I$, that is $ L( A;D) > J / I \Omega(A) = \Lambda $ where your critical likelihood ratio $\Lambda$ is denied by the context.
Classical Statisticians were the first to conclude that there must be some $\Lambda$ such that you will guess $A$ if
$ L( A;D)> \Lambda $ and guess $\bar{A}$ if  $ L( A;D)< \Lambda $. By and large, classical statisticians say the choice of $\Lambda$
is an entirely subjective one which no one but you can make (e.g. Lehman, 1959, p.62). 
Bayesians agree; $\Lambda$ is inversely proportional to your current
odds for $A$, an aspect of your personal opinion.
The classical statisticians, however, have overlooked a great simplification, namely that your critical $\Lambda$ will not depend on the size or
structure of the experiment and will be proportional to $J/I$. 
As Savage (1962) puts it: the subjectivist's position is more objective than the objectivist's, for the subjectivist finds the range of coherent
or reasonable preference patterns much narrower than the objectivist thought it to be. How confusing and dangerous big words are (p.67)!}
\end{quote}

Given this discussion, we build on the idea that a hypothesis test can be constructed by focusing on the distribution of the test statistic, denoted by $t$, under the alternative hypothesis.  Bayes factors can then  be calculated once the researcher is willing to assess a prior predictive interval for the $t$ statistic under the alternative.  In most experimental situations, this appears to be the most realistic way of assessing {\em a priori} information.  For related discussion, see Berger and Sellke (1987) and Berger (2003) who pose the question of whether Fisher, Jeffreys and Neyman could have agreed on testing  and provide illuminating examples illustrating the differences (see Etz and Wagenmakers, 2017).

The rest of our paper is outlined as follows.  Section 2 provides a framework for the differences between Classical and Bayesian hypothesis testing.  Section 3 uses a probabilistic interval assessment for the test statistic distribution under the alternative to assess a Bayes factor.  Jeffrey's (1957, 1961) Cauchy prior and the Bartlett-Lindley paradox (
Lindley, 1957, and Bartlett, 1957) are discussed in this context.
Extensions to regression and $R^2$, $\chi^2$ and $F$ tests (see Connely, 1991, and Johnson, 2005, 2008) are also provided.  Section 4 concludes with further discussion and with directions for future research.

\section{Bayesian vs Classical Hypothesis Testing}
Suppose that you wish to test a sharp null hypothesis $ H_0 : \theta=0$ against a non-sharp composite alternative
$ H_1 : \theta \neq 0 $.   We leave open the possibility that $H_0$ and $H_1$ could represent models and the researcher wishes to perform model selection.  A classical test procedure uses the sampling distribution, denoted by $p(\hat{\theta}|\theta)$, of a test statistic ${\hat \theta}$, given the parameter $\theta$.  A critical value, $c$, is used to provide a test procedure of the form
$$
{\rm Reject} \;  H_0 \; {\rm if} \; | \hat{\theta} | > c.
$$
There are two types of errors that can arise. Either the hypothesis maybe rejected even though it is true (a Type I error)
or it maybe accepted even though it is false (Type II). Typically, the critical value $c$ is chosen so as to make the probability of a type I error, $\alpha$, to
be of fixed size.  We write $ \alpha(c) = 1 - \int_{-c}^c p( \hat{\theta} | \theta ) d \hat{\theta} $. 

Bayes factor, denoted by BF, which is simply a density ratio (as opposed to a tail probability) is defined by a likelihood ratio
$$
BF = \frac{p(\hat{\theta} |H_0 )}{p(\hat{\theta}|H_1 )}.
$$
Here $p({\hat \theta}|H_0)=\int p({\hat \theta}|\theta,H_0)p(\theta|H_0)d\theta$ is a marginal distribution
of the test statistic and $p(\theta|H_0)$ an {\em a priori} distribution on the parameter.  For a simple hypothesis, 
$(\theta|H_0) \sim \delta_{\theta_0}$ is a Dirac measure at the null value.  The difficulty comes in specifying $p(\theta|H_1)$, the prior under the alternative.  A Bayesian Hypothesis Test can then be constructed in conjunction with the {\em a priori} odds ratio $p(H_0)/p(H_1)$, to calculate a posterior odds ratio, via Bayes rule,
$$
\frac{p(H_0|\hat{\theta} )}{p(H_1|\hat{\theta} )} = \frac{p(\hat{\theta} |H_0 )}{p(\hat{\theta}|H_1 )}\frac{p(H_0 )}{p(H_1 )}.
$$
As $ H_0 : \theta=\theta_0$ and $ H_1 : \theta \neq \theta_0 $, the Bayes factor calculates $p({\hat \theta}|\theta=\theta_0)/p({\hat \theta}|\theta \neq \theta_0)$.  
We will focus on the test statistic distribution under the alternative hypothesis, namely $p({\hat \theta}|\theta \neq \theta_0)$.
See also Held and Ott (2018) for additional discussion on $p$-values and Bayes factors.

\section{A Default Bayes Factor}

Our approach is best illustrated with the usual $t$-ratio test statistic in a normal means problem.
As ELS illustrate, the central question that a Bayesian must \emph{a priori} address is the question:
\begin{quote}
\emph{How often will $t$ lie in that interval when the null hypothesis is false?}
\end{quote}
To do this we need an assessment to the distribution of the $t$-ratio test statistic under the alternative.
As ELS further observe:
\begin{quote}
\it{This argument, like almost all the
following discussion of null hypothesis
testing, hinges on assumptions about
the prior distribution under the alternative
hypothesis. The classical statistician
usually neglects that distribution
in fact, denies its existence.
He considers how unlikely a $t$ as far
from $0$ as $1.96$ is if the null hypothesis
is true, but he does not consider that
a t as close to $0$ as $1.96$ may be even
less likely if the null hypothesis is
false.}
\end{quote}

First, we calculate prior predictive distribution of the test statistic under the alternative and then show how such assessment
can lead to a default Bayes factor.

\subsection{Predictive distribution, $Pr(T=t|H_1)$}
A simple default approach to quantifying {\em a priori} opinion is to assess a hyperparameter, denoted by $A$, such that the following probability
statements hold true:
\begin{align*}
Pr \left ( - 1.96 \sqrt{A} < T < 1.96 \sqrt{A} | H_1 \right ) & = 0.95\\
Pr \left ( - 1.96 < T < 1.96 | H_0 \right ) & = 0.95.
\end{align*}
Under the null, $H_0$, both the Bayesian and Classicist agree that $A=1$.  All that is needed to complete the specification is the assessment of $A$.

In the normal mean testing problem we have an iid sample $( y_i |\theta) \sim N(\theta, \sigma^2)$, for $i=1,\ldots,n$,  with $\sigma^2$ known and  $n{\bar y}=\sum_{i=1}^n y_i$.  Under the null, $H_0: \theta=0$, the distribution of $T = \sqrt{n} \bar{y} / \sigma $, is the standard normal distribution, namely $T \sim N(0,1)$.  The distribution of $T$ under the alternative, $H_1: \theta \neq 0$, is a mixture distribution
$$
p(T=t|H_1 ) = \int_\Theta p( T=t|\theta)p(\theta|H_1)d \theta,
$$
where $p(\theta|H_1)$ denotes the prior distribution of the parameter under the alternative.
Under a normal sampling scheme, this is a location mixture of normals
$$
p( T=t|H_1) = \int_{-\infty}^\infty P(T=t|H_1,\theta)p(\theta|H_1)d \theta
$$
where $T|H_1,\theta$ is normal with mean $\sqrt{n}\theta/\sigma$ and variance one; or $T=\sqrt{n}\theta/\sigma + \varepsilon$, where $\varepsilon \sim N(0,1)$.

Under a normal prior, $\theta \sim N( 0 , \tau^2)$, the distribution $p(T=t|H_1)$ can be calculated in closed form as
 $T \sim N(0,A) $ where $A=1+n\tau^2/\sigma^2$.  Hence an assessment of $A$ will depend on the design (through $n$) and the relative ratio of measurement errors (through $\tau^2/\sigma^2$).  The gain in simplicity of the Bayes test is off-set by the difficulty in assessing $A$.

The Bayes factor is then simply the ratio of two normal ordinates
$$
B = \frac{ \phi(t)}{ \phi( t/\sqrt{A}) } = \sqrt{A} \exp\left\{-\frac{1}{2} t^2(1-A^{-1})\right\} \;.
$$
The factor $A$ is often interpreted as the Occam factor (Berger and Jefferys, 1992,  
Jefferys and Berger, 1992, Good, 1992).  See Hartigan (2003) for a discussion of default Akaike-Jeffreys priors and model selection.

Our approach requires the researcher to ``calibrate'' $A$ ahead of time.  One simple approach is to perform a \emph{what if} analysis and assess what posterior
odds we would believe \emph{if} we saw $t=0$. 
This assessment directly gives the quantity $\sqrt{A}$.

\paragraph{Dickey-Savage.}
The Bayes factor $BF$ for testing $ H_0 $ versus $ H_1 $ can be calculated using the 
Dickey-Savage density ratio. This relates the posterior model probability $p( \theta= \theta_0 | y)$ to
the marginal likelihood ratio via Bayes rule
$$
\frac{Pr( \theta= \theta_0 | y) }{Pr( \theta= \theta_0 ) } = \frac{ p( y | \theta= \theta_0 ) }{ p( y ) }.
$$

\paragraph{Bayes Factor Bounds.}
Let ${\hat \theta}_{MLE}$ denote the maximum likelihood estimate, then
$$ 
p( y | H_1) = \int p( y | \theta ) p( \theta | H_1 ) d \theta \leq p( y | \hat{\theta}_{MLE} ).
$$ 
This implies that, for $H_0: \theta=0$,
$$
BF \geq \frac{p(T=t|\theta=0)}{p(T=t|\hat{\theta})}.
$$
In a normal means testing context, this leads to a bound,
$$
\frac{ p( y | H_0 ) }{ p( y | H_1 )}  \geq \exp\{-0.5(1.96^2-0^2)\} = 0.146\;.
$$
Under an {\em a priori} $1/2$-$1/2$ weight on either side of zero, the bound increases to $0.292$. 
Hence, \emph{at least 30\% of the hypotheses that the classical approach rejects are true in the Bayesian world}.
Amongst the experiments with $p$-values of $0.05$ at least $30$\% will actually turn out to be true!
Put another way, the probability of rejecting the null \emph{conditional} on the observed $p$-value of $0.05$ is at least $30$\%.
You are throwing away good null hypothesis and claiming you have found effects!
In terms of posterior probabilities, with $p(H_0)=p(H_1)$, we have a bound
$$
Pr( H_0 | y ) = \left[1 + \frac{ p( y | H_1 ) }{ p( y| H_0 )} \frac{ Pr( H_1 ) }{ Pr( H_0 )} \right]^{-1} \geq 0.128\;.
$$
Hence, there is at least $12.8$ percent chance that the null is still true even in the one-sided version of the problem!   
Clearly at odds with a $p$-value of $5$ percent.

One of the key issues, as discussed by ELS, is that  the classicist approach is based on an observed $p$-value is not a probability in any real sense. The observed $t$-value is a realization of a statistic
that happens to be $N(0,1)$ under the null hypothesis. Suppose that we observe $t=1.96$.
Then the \emph{maximal evidence} against the null hypothesis which corresponds to $t=0$ will be achieved by evaluating
the likelihood ratio at the observed $t$ ratio, which is distributed $N(0,1)$. 

\subsection{Normal means Bayes factors}
We have the following set-up for the normal means case
(see Berger and Delampaday, 1989, for the full details): Let 
${\bar y}|\theta \sim N(\theta,\sigma^2/n)$, where $\sigma^2$ is known and let
$t =\sqrt{n}({\bar y}-\theta_0)/\sigma$ the t-ratio test statistic when 
testing the null hypothesis $H_0: \theta = \theta_0$ against the alternative hypothesis $H_0: \theta \neq \theta_0$.  Also, 
let $ \rho = \sigma / ( \sqrt{n} \tau ) $
and $ \eta = ( \theta_0 - \mu )/\tau $, derived from a normal prior in the alternative $ \theta
\sim N (\mu, \tau^2 ) $. Usually, we take a symmetric prior and set $\mu = \theta_0 $, such that
$ \eta = 0$ and the Bayes factor simplifies to
$$
BF = \sqrt{ 1 + \rho^{-2} } \exp \left ( - \frac{1}{2 (1 + \rho^2 )} t^2 \right ).
$$
We can use the Dickey-Savage density ratio as follows to derive the above Bayes factor:
\begin{align*}
p(\theta_0 | {\bar y}) & = \frac{1}{\sqrt{2 \pi} \tau \sqrt{ 1 + \rho^{-2} }}
 \exp \left ( - \frac{1}{2 (1 + \rho^2 )} t^2 \right ) \\
p(\theta_0) & = \frac{1}{\sqrt{2 \pi} \tau}
\end{align*}
The posterior distribution under the alternative is
$$
( \theta | y ) \sim \mathcal{N} \left (  \left ( \frac{n}{\sigma^2} + \frac{1}{\tau^2} \right )^{-1}  
 \left ( \frac{n \bar{y}}{\sigma^2} + \frac{\theta_0}{\tau^2} \right )
 , \left ( \frac{n}{\sigma^2} + \frac{1}{\tau^2} \right )^{-1} \right )
$$
with quantities
$$
t^2 = \frac{n ( \bar{y} - \theta_0 )^2 }{\sigma^2} \;
{\rm and} \; \left ( \frac{n}{\sigma^2} + \frac{1}{\tau^2} \right )^{-1} = \tau^2 ( 1 + \rho^{-2} )^{-1}.
$$
The posterior mean $E( \theta | y)$ can be written as  
$$
\theta_0 + \left ( \frac{T}{\sigma^2} + \frac{1}{\tau^2} \right )^{-1} \frac{T ( \bar{y} - \theta_0 ) }{\sigma^2}.
$$
Substituting into the ratio of ordinates $ p( \theta= \theta_0 | y) /p( \theta= \theta_0) $ gives the result.

In the case where $ \tau$ is moderate to large, this is approximately
$$
BF = \frac{ \sqrt{n} \tau }{\sigma }  \exp \left ( - \frac{1}{2} t^2 \right ).
$$
Clearly, the prior variance $\tau$ has a dramatic effect on the answer. 
First, we can see that the ``noninformative'' prior $ \tau^2 \rightarrow \infty $ makes
little sense (Lindley, 1957, Bartlett, 1957).  For instance, when $\sigma=\tau$ and $t=2.567$ (a $p$-value of 0.01), 
then the Bayes factor equals
$0.16$, $1.15$, $3.62$ and $36.23$ for $n$ equal to $10$, $100$, $1000$ and $1000000$, respectively 
(see Section 3.3 for more details about the Bartlett-Lindley Paradox).
Secondly, the large effect is primarily
due to the thinness of the normal prior in the tails. Jeffreys (1961) then proposed
the use of a Cauchy $(0,\sigma^2)$ prior (see Section 3.4 for further details).

\subsection{Bartlett-Lindley Paradox}
See Lindley (1957) and Bartlett (1957) for the full details.
The Barlett-Lindley paradox occurs when you let $ \tau^2 \rightarrow \infty$. This has the ``appropriate'' behaviour at the origin of 
flattening out the marginal distribution of $T$. So when comparing equal length intervals
$Pr( a < T < b )$ and $Pr( c < T < d ) $, where $ a-b=c-d$, one would get approximately a Bayes factor of one.

The so-called paradox arises when the Bayes factor places all its weight on the alternative hypothesis $ H_1$.
Thought of via the marginal predictive of $T$ this is not surprising. As $ \tau^2 \rightarrow \infty$ implies
$ A \rightarrow \infty$, and your belief {\em a priori} that you expect an incredibly large value of $T$ values under the alternative.
Now, when you actually observe $ 1.96<T<2.56$ this is unlikely under the null approximately $2$\%, but nowhere near as likely under the alternative. The Bayes factor correctly identifies the null as having the most posterior mass.

\subsection{Cauchy Prior}
Jeffreys (1961) proposed a Cauchy 
(centered at $\theta_0$ and scale $1$) to allow for fat-tails
whilst simultaneously avoiding having to specify a scale to the normal prior.
Using the asymptotic, large $n$, form of the posterior $(\sqrt{n}/\sqrt{2 \pi}\sigma)\exp\{-0.5 t^2\}$ for the usual $t$-ratio test statistic and the fact that the prior
density ordinate from the Cauchy prior is $p(\theta_0 ) = 1/(\pi\sigma)$, the Bayes Factor is
$$
BF = \frac{(\sqrt{n}/\sqrt{2 \pi}\sigma)\exp\{-0.5 t^2\}}{1/(\pi\sigma)} = \sqrt{0.5\pi n}\exp\{-0.5 t^2\}.
$$
We have the interval probability
$$
Pr( - 1.96 \sqrt{A}< T < - 1.96 \sqrt{A} | H_1)  \approx 0.95,
$$
for $A \approx 40$, when $n=1$ and $\sigma^2=1$.  Exact answer given by cdf of hypergeometric Beta You can also see this in the Bayes factor approximations.   
Therefore, very different from letting $A \rightarrow \infty$, in a normal prior.

\subsection{Coin tossing: $p$-values and Bayes}

Suppose that you routinely reject two-sided hypotheses at a fixed level of significance, say $\alpha =0.05$.
Furthermore, suppose that half the experiments under the null are actually true, i.e. $Pr(H_0)=Pr(H_1)=0.5$.
The experiment will provide data, $y$, here we standardize the mean effect and obtain a $t$-ratio.

\paragraph{Example: Coin Tossing (ELS).}
Let us start with a coin tossing experiment where you
want to determine whether the coin is ``fair", $H_0: Pr(Head) = Pr(Tail)$, or 
the coin is not fair, $H_1: Pr(Head) \neq Pr(Tail)$.
ELS discuss at length the following four experiments where, in each case, the test statistics is $t=1.96$.
We reproduce below of their Table 1.

\begin{table}[H]
\begin{center}
\begin{tabular}{rrrrr}\hline
Expt & 1 & 2 & 3 & 4 \\  \hline
$n$ & 50 & 100 & 400 & 10,000 \\
$r$ & 32 & 60 & 220 & 5,098 \\
$BF$ & 0.8 & 1.1 & 2.2 & 11.7 \\ \hline
\end{tabular}
\end{center}
\caption{The quantities $n$ and $r$ are, respectively, number of tosses of the coin and the number of heads that 
barely leads to rejection of the null hypothesis, $H_0: Pr(Head) = Pr(Tail)$, by a classical two-tailed test at the 5 percent level.}
\end{table}
For $n$ coin tosses and $r$ heads, the Bayes factor,
$$
BF = \left(\frac12\right)^n/\int_0^1 \theta^r(1-\theta)^{n-r}p(\theta|H_1)d\theta,
$$
which grows to infinity and so there is overwhelming evidence in favor of  $H_0: Pr(Head) = Pr(Tail)$.  
This is a clear illustration of Lindley's paradox.

There are a number of ways of assessing the odds. One is to use a uniform prior. Another useful approach which gives a lower bound
is to use the \emph{maximally informative} prior which puts all its mass on the parameter value at the mle, $\hat{\theta}=r/n$.
For example, in the $r=60$ versus $n=100$ example, we have $ \hat{\theta} = 0.6 $. Then we have 
$ p( y | H_1 ) \leq p( y | \hat{\theta} ) $ and for the odds ratio
$$
\frac{ p( y | H_0 ) }{ p( y | H_1 )} \geq \frac{ p( y | \theta = \theta_0 ) }{ p( y| \hat{\theta} ) }.
$$
For example, with $ n=100$ and $r=60$, we have
$$
\frac{ p( y | H_0 ) }{ p( y | H_1 )} \geq \frac{0.5^{100} }{ 0.6^{60} 0.4^{40} } = 0.134.
$$
In terms of probabilities, if we start with a $50/50$ prior on the null, then the posterior probability of the null is at least $0.118$:
$$
Pr(H_0|y) = \left(1 + \frac{ p( y| H_1 ) }{ p( y | H_0 )} \frac{Pr(H_1)}{Pr(H_0)} \right)^{-1} \geq 0.118.
$$

\subsection{Regression}
A number of authors have provided extensions to traditional classical tests, for example 
Johnson (2008) shows that $R^2$, deviance, $t$ and $F$ can all be interpreted as Bayes factors.
See also Gelman {\em et al} (2008) for weakly informative default priors for logistic regression models.

In the case of nested models, Connelly (1991) proposes the use of
$$
BF = n^{ - \frac{d}{2} } \left ( 1 + \frac{d}{n-k} F \right )^{ \frac{n}{2} }
$$
where $ F$ is the usual $F$-statistic, $k$ is the number of parameters in the larger model
and $d$ is the difference in dimensionality between the two models.
In the non-nested case, first it helps to nest them if you can, otherwise MCMC comparisons.

Zellner and Siow (1980) extend this to the Cauchy prior case, see Connelly (1991).
Essentially, introduces a constant out-front that depends on the prior ordinate $p(\theta=\theta_0)$.
See Efron and Gous (2001) for additional discussion of model selection in the Fisher and Jeffreys approaches.
Additionally, Polson and Roberts (1994) and  Lopes and West (2004) study model selection in diffusion processes and
factor analysis, respectively.  Scott and Berger (2010) compare Bayes and empirical-Bayes in the variable selection context.

\section{Discussion}

The goal of our paper was to revisit ELS.  
There are a number of important take-aways from comparing the Bayesian paradigm to frequentist ones.  
Jeffreys (1937) provided the foundation for Bayes factors (see Kass and Raftery, 1995, for a review).
Berkson (1938) was one of the first authors to point out problems with p-values.

The Bayesian viewpoint is clear: you have to \emph{condition} on what you see.
You also have to make probability assessments about competing hypotheses.
The \emph{observed} $y$ can be highly unlikely under \emph{both} scenarios! 
It is the relative oods that is important.
The $p$-value under both hypotheses are then very small, but the Bayes posterior probability is based on the \emph{relative} odds of observing the data
plus the prior, that is $p(y| H_0) $ and $p(y|H_1)$ can both be small, but its $ p( y| H_0)/p(y|H_1) $ that counts together with the prior $p(H_0)/p(H_1)$. 
Lindley's paradox shows that a Bayes test has an extra factor of $\sqrt{n} $ which will asymptotically
favor the null and thus lead to asymptotic differences between the two approaches.  
There is only a practical problem when $ 2 < t < 4 $ -- but this is typically the most interesting case!

Jeffreys (1961), page 385, said that ``{\em what the use of P implies \ldots is that a hypothesis that may be true may be rejected because it has not predicted observable results that have not occurred.  This seems a remarkable procedure.}''

We conclude with two quotes on what is wrong with classical p-values 
with some modern day observations from two Bayesian statisticians.

\paragraph{Jim Berger:} {\em 
$p$-values are typically much smaller than actual error probabilities
$p$-values do not properly seem to reflect the evidience in the data. For instance, suppose one pre-selected $\alpha=0.001$.
This then is the error one must report whether $p=0.001$ or $p=0.0001$, in spite of the fact that the latter would seem to provide much stronger evidence against
the null hypothesis.}

\paragraph{Bill Jefferys:} {\em 
The Lindley paradox goes further. It says, assign priors however you wish. You don't get to change them.
Then take data and take data and take data ... There will be times when the classical test will reject with probability
$(1-\alpha)$ where you choose $\alpha$ very small in advance, and at the same time the classical test will reject at a significance level $\alpha$.
This will not happen, regardless of priors, for the Bayesian test. The essence of the Lindley paradox is that ``sampling to a foregone conclusion''
happens in the frequentist world, but not in the Bayesian world.}

As we pointed out at the outset, hypothesis testing is still a central issue in modern-day statistics and machine learning, in particular, its relationship with high dimensional model selection.  Finding default regularization procedures in high dimensional settings is still an attractive area of research.

\section*{References}

\begin{description}
\item Bartlett, M. S. (1957) Comment on ``A Statistical Paradox'' by D.V. Lindley.
\textit{Biometrika}, 44, 533-534.
\medskip

\item Berger, J. O. (2003) Could Fisher, Jeffreys and Neyman have agreed on testing?
\textit{Statistical Science}, 18, 1-32.
\medskip

\item Berger, J. O. and Delampady, M. (1987) Testing Precise Hypothesis (with Discussion).
\textit{Statistical Science}, 2, 317-335.
\medskip

\item Berger, J. O. and Jefferys, W. H. (1992) The application of robust Bayesian analysis to hypothesis testing and Occam's Razor.  {\em Journal of the Italian Statistical Society}, 1, 17-32.

\item Berger, J. O. and Sellke, T. (1987) Testing of a point null hypothesis: the irreconcilability of significance levels and evidence (with Discussion). {\em Journal of the American Statistical Association}, 82, 112-139.

\item Berkson, J. (1938) Some difficulties of interpretation encountered in the application of the chi-square test.  {\em Journal of the American Statistical Association},
33, 526-542.

\item Connelly, R. A. (1991) A Posterior Odds Analysis of the Weekend Effect.
\textit{Journal of Econometrics}, 49, 51-104.
\medskip

\item Edwards, W., Lindman, H.  and Savage, L. J. (1963)  Bayesian Statistical Inference for Psychological Research.
\textit{Psychological Review}, 70(3), 193-242.
\medskip

\item Efron, B. and Gous, A. (2001) 
Scales of Evidence for Model Selection: Fisher versus Jeffreys. {\em Model selection}, 208-246, Institute of Mathematical Statistics, Beachwood.
\medskip

\item Etz, A. and Wagenmakers, E.-J. (2017) J. B. S. Haldane's Contribution to the Bayes Factor Hypothesis Test.
{\em Statistical Science}, 32(2), 313-329.

\item Gelman, A., Jakulin, A., Pitau, M. G. and Su, Y.-S. (2008) A Weakly Informative Default Prior Distribution for Logistic and Other Regression Models. 
{\em The Annals of Applied Statistics}, 2(4), 1360-1383.
\medskip

\item Good, I. J. (1992) The Bayes/Non-Bayes Compromise: A Brief Review.
{\em Journal of the American Statistical Association}, 87, 597-606.
\medskip

\item Held, L. and Ott, M. (2018) On $p$-Values and Bayes Factors.  
{\em Annual Review of Statistics and Its Application}, 5, 393-419.
\medskip

\item Jefferys, W. H.  and Berger, J. O. (1992) Ockham's Razor and Bayesian Analysis.
\textit{American Scientist}, 80, 64-72.

\item Jeffreys, H. (1957) {\em Scientific inference}. Cambridge: Cambridge University Press.

\item Jeffreys, H. (1961) \textit{Theory of Probability}. London: Oxford University Press.
\medskip

\item Johnson, V. (2005) Bayes factors based on test statistics.
{\em Journal of the Royal Statistical Society, Series B}, 67, 689-701.

\item Johnson, V. (2008) Properties of Bayes Factors Based on Test Statistics.
{\em Scandinavian Journal of Statistics}, 35, 354-368.

\item Kass, R. E. and Raftery, A. E. (1995) Bayes factors.  {\em Journal of the American Statistical Association}, 90, 773-395.

\item Lehmann, E. L. (1959) {\em Testing statistical hypotheses}, New York: Wiley.

\item Lindley, D. V. (1957). A Statistical Paradox. \textit{Biometrika}, 44, 187-192.
\medskip

\item Lopes, H. F. and West, M. (2004) Bayesian model assessment in factor analysis.
{\em Statistica Sinica}, 14, 41-67.
\medskip

\item Polson, N. G. and Roberts, G. O. (1994) Bayes Factors for discrete observations from 
diffusion processes.  {\em Biometrika}, 81(1), 11-26.
\medskip

\item Savage, L. J. (1962) Subjective probability and statistical practice. In L. J. Savage et al.,
{\em The foundations of statistical inference: A discussion}.  New York: Wiley.

\item Scott, J. G. and Berger, J. O. (2010) Bayes and empirical-Bayes multiplicity adjustment in the variable-selection problem.  
{\em Annals of Statistics}, 38(5), 2587-2619.
\medskip

\item Zellner, A. and Siow, A. (1979) Posterior Odds Ratio for Selected Regression Hypotheses.
In: \textit{Bayesian Statistics}, Proceedings of the First International Meeting (J. M. Bernardo, M. H. De Groot, D. V. Lindley and A. F. M. Smith, eds), pp. 585-603. Valencia: University Press.

\end{description}
 
\end{document}